\documentclass[10pt]{amsart}

\usepackage{mathabx}
\input bookman.sty
\boldmath

\usepackage{comment}
\excludecomment{excl}

\usepackage{xcolor}

\usepackage{helvet}

\usepackage{graphics}
\usepackage{amssymb}
\usepackage{amsxtra}
\usepackage{amsmath}
\usepackage{mathrsfs}

\usepackage[arrow, matrix]{xy}
\xyoption{frame}

\theoremstyle{plain}

\newtheorem{theo}{Theorem}[section]
\newtheorem{lemm}[theo]{Lemma}
\newtheorem{prop}[theo]{Proposition}
\newtheorem{coro}[theo]{Corollary}

\theoremstyle{definition}

\newtheorem{defi}[theo]{Definition}
\newtheorem{rema}[theo]{Remark}

\newfont{\rmm}{cmr10 scaled 1000}

\newfont{\itt}{cmsl10 scaled 1000}

\newfont{\rM}{cmr10 scaled 1700}

\newcounter{lemma}[section]

\newcounter{tempcounter}

\newcommand{\lb}{\label}

\newcommand{\rrf}[1]{(\ref{#1})}

\renewcommand{\a}{\alpha}
\renewcommand{\b}{\beta}
\newcommand{\g}{\gamma}
\renewcommand{\d}{\delta}

\renewcommand{\l}{\lambda}

\newcommand{\m}{\mu}

\renewcommand{\r}{\rho}

\newcommand{\s}{\sigma}

\renewcommand{\o}{\omega}

\newcommand{\G}{\Gamma}

\renewcommand{\L}{\Lambda}

\renewcommand{\O}{\Omega}

\newcommand{\cc}{{\mathbb{C}}}

\newcommand{\rr}{{\mathbb{R}}}

\newcommand{\ttt}{{\mathbb{T}}}

\newcommand{\zz}{{\mathbb{Z}}}

\newcommand{\MMM}{{\mathbf{M}}}

\newcommand{\CCCC}{{\mathscr{C}} }
\newcommand{\DDDD}{{\mathscr{D}} }
\newcommand{\EEEE}{{\mathscr{E}} }

\newcommand{\MMMM}{{\mathscr{M}}}

\newcommand{\id}{\text{id}}

\newcommand{\Ext}{\text{\rm Ext}}

\newcommand{\bere}{\begin{rema}}
\newcommand{\bede}{\begin{defi}}

\renewcommand{\beth}{\begin{theo}}
\newcommand{\bele}{\begin{lemm}}
\newcommand{\bepr}{\begin{prop}}
\newcommand{\beeq}{\begin{equation}}
\newcommand{\bega}{\begin{gather}}
\newcommand{\begaa}{\begin{gather*}}
\newcommand{\been}{\begin{enumerate}}

\newcommand{\bedee}{\begin{defii}}
\newcommand{\bethh}{\begin{theoo}}
\newcommand{\belee}{\begin{lemmm}}
\newcommand{\beprr}{\begin{propp}}

\newcommand{\beco}{\begin{coro}}

\newcommand{\beal}{\begin{aligned}}

\newcommand{\enre}{\end{rema}}

\newcommand{\enco}{\end{coro}}
\newcommand{\enpr}{\end{prop}}
\newcommand{\enth}{\end{theo}}
\newcommand{\enle}{\end{lemm}}
\newcommand{\enen}{\end{enumerate}}
\newcommand{\enga}{\end{gather}}
\newcommand{\engaa}{\end{gather*}}
\newcommand{\eneq}{\end{equation}}
\newcommand{\enal}{\end{aligned}}

\newcommand{\bq}{\begin{equation}}
\newcommand{\bqq}{\begin{equation*}}

\renewcommand{\geq}{\geqslant}

\newcommand{\wi}{\widetilde}

\newcommand{\ove}{\overline}

\newcommand{\wh}{\widehat}

\newcommand{\sbs}{\subset}

\newcommand{\sut}{~such~that~}

\newcommand{\ho}{homomorphism}

\newcommand{\ma}{manifold}

\newcommand{\Prf}{{\it Proof.\quad}}

\newcommand{\smo}{C^{\infty}}

\newcommand{\chart}{\Phi_p:U_p\to B^n(0,r_p)}
\newcommand{\atlas}{\{\Phi_p:U_p\to B^n(0,r_p)\}_{p\in S(f)}}

\newcommand{\pr}{\partial}

\newcommand{\qs}{\hfill\square}

\newcommand{\arrh}[3]
{
\xymatrix{
{#1} \ar[r]^<<<<{#2}  &{#3}
}
}

\newcommand{\arrr}[1]
{\arrh {}{#1}{}}

\newcommand{\arrto}
{\xymatrix{{} \ar@{|-{>}}[r]  & {} } }

\newcommand{\arrinto}
{\xymatrix{{} \ar@{^{(}->}[r]  & {} } }

\newcommand{\ses}{spectral sequence}

\begin{document}

\title
{Massey products in mapping tori}
\author{Andrei Pajitnov}
\address{Laboratoire Math\'ematiques Jean Leray 
UMR 6629,
Universit\'e de Nantes,
Facult\'e des Sciences,
2, rue de la Houssini\`ere,
44072, Nantes, Cedex}                    
\email{andrei.pajitnov@univ-nantes.fr}

\thanks{} 
\begin{abstract}

Let $\phi: M\to M$ be a diffeomorphism of a $\smo$ compact connected \ma~,
and $X$ its mapping torus. There is a natural fibration
$p:X\to S^1$, denote by $\xi\in H^1(X, \zz)$
the corresponding cohomology class. Let $\l\in \cc^*$.

Consider the endomorphism 
$\phi_k^*$
induced by $\phi$ in the cohomology of $M$ of degree $k$,
and denote by $J_k(\l)$ 
the maximal size of its Jordan block of eigenvalue $\l$.
Define a representation
$\r_\l : \pi_1(X)\to\cc^*; \ \r_\l(g)=\l^{p_*(g)}$;
let $H^*(X,\r_\l)$ be the corresponding twisted cohomology of $X$.
We prove that $J_k(\l)$ is equal to the maximal length 
of a non-zero Massey product of the form
$\langle \xi, \ldots , \xi, a\rangle$
where $a\in H^k(X,\r_\l)$
(here the {\it length } means the number of entries of $\xi$).

In particular, if $X$ is a strongly formal space (e.g. a K\"ahler manifold)
then all the Jordan blocks of 
$\phi_k^*$
are of size 1. 
If $X$ is a formal space, then all the Jordan blocks of eigenvalue 1 
are of size 1. This leads to a simple construction
of formal but not strongly formal mapping tori.

The proof of the main theorem is based on the fact that 
the Massey products of the above form can be identified with differentials
in a Massey spectral sequence, which  in turn 
can be explicitly computed in terms of the 
Jordan normal form of $\phi^*$.

\end{abstract}
\maketitle
\tableofcontents

\section{Introduction}
\label{s:intro}


The relation between non-vanishing Massey products of length 2 and 
the Jordan blocks of size greater than 1 was discovered 
in the work of  M. Fern\'{a}ndez, A. Gray, J. Morgan,  \cite{FGM},
where it was used to prove that certain mapping tori 
do not admit a structure of a K\"ahler manifolds.
In the work of G. Bazzoni,  M. Fern\'{a}ndez, V. Mu\~{n}oz \cite{BFM}
it was proved that the existence of Jordan blocks of size 2  
implies the existence of a non-zero triple Massey product of the form
$\langle  \xi, \xi, a\ \rangle$.

The main theorem of the present paper provides systematic treatment 
of these phenomena,
relating the length of non-zero Massey products to the size of 
Jordan blocks.
The both numbers turn out to be equal 
to the number of the sheet where the formal deformation
spectral sequence degenerates.

Another approach to the relation between the size of Jordan blocks and formality 
was developed by S. Papadima and A. Suciu \cite{PS}, \cite{PS2}.
They prove in particular that if the monodromy \ho~
has Jordan blocks of size greater than 1, then the 
fundamental group of the mapping torus is not a formal group.

\subsection{Overview of the article}
\label{su:overv}
The proof of the main theorem is based on the 
techniques developed in \cite{KP}.
We begin by an overview of this paper in Section \ref{s:sp-seq}.
The main theorem of the paper is stated and proved in 
Section \ref{s:m_th}.
In Section \ref{s:gener}
we present a generalization of the main theorem to the case
of spaces $X$ endowed with a non-zero
cohomology class $\xi\in H^1(X,\zz)$.

\section{Formal deformations and Massey spectral sequences.}
\lb{s:sp-seq}

Let $X$ be a connected manifold, and $\xi\in H^1(X, \cc)$ a non-zero cohomology class,
and $\l\in \cc^*$.
There is a spectral sequence starting with $H^*(X,\r_\l)$ 
and converging to the cohomology $H^*(X,\r_{\l'})$ where $\l'$ is a generic complex number.
There are different versions of this spectral sequence in literature,
see \cite{F1}, \cite{NovM}, \cite{PaM}, \cite{F2}, \cite{KP}.
We will recall here the versions described in \cite{KP},
refering to this article for details and proofs. 

\subsection{Massey spectral sequences}
\label{su:mas-ss}

In this subsection we describe a spectral sequence with the differential 
defined in terms of special Massey products. 

Pick $\a\in\cc$ such that $e^\a=\l$. 
The cohomology  $H^*(X,\r_\l)$ can be computed  from 
the twisted DeRham complex 
$\tilde\Omega^*(X) = \Omega^*(X, \tilde d)$ where 
$\tilde d(\o)=d\o+\a\xi\wedge \o$.
Let $a\in  H^*(X,\r_\l)$. 
An {\it $r$-chain starting from $a$} is a sequence of differential forms 
$\o_1, \ldots, \o_r\in \Omega^*(X)$
such that 
$$
 d\o_1=0,\ \  [\o_1]=a,\ \   d\o_2=\xi\wedge\o_1,\ \ldots,\   d\o_r=\xi\wedge\o_{r-1}.
$$
For an $r$-chain $C$ put $\pr C=\xi\wedge\o_r$; this is a cocycle in 
$\tilde\Omega^*(X)$
.
Denote by 
$MZ^m_{(r)}$
 the subspace of all 
$a\in   H^*(\tilde\Omega^*(X))$
 such that there exists an $r$-chain
starting from $a$. 


Denote by 
$MB^m_{(r)}$
the subspace of all 
$\beta\in   H^*(\tilde\Omega^*(X))$ such that there 
exists an $(r-1)$-chain
$C=(\o_1, \ldots , \o_{r-1})$ with $\xi\wedge\o_r$ belonging to $\beta$. 
It is clear that 
$MB^m_{(i)}
\sbs
MZ^m_{(j)}$
for every $i, j$.
Put
$$
MH^m_{(r)}
=
MZ^m_{(r)} \Big/
MB^m_{(r)}.
$$
In the next definition we  omit the upper indices and write
 $MH_{(r)},MZ_{(r)} $ etc.  in order to simplify the notation.
\begin{defi}
 Let $a\in  H^*(\tilde\Omega^*(X))$, and $r\geq 1$.
We say that the {\it $r$-fold Massey product
$\langle \xi, \ldots, \xi,a\rangle$
is defined, if $a\in MZ_{(r)}$.}
In this case choose any $r$-chain 
$(\o_1, \ldots , \o_{r})$ 
starting from $a$. The cohomology class of $\pr C=\xi\wedge\o_{r}$
is in $MZ_{(r)}$
and  is well defined 
modulo $MB_{(r)}$.
The image of $\pr C$ in $MZ_{(r)}/MB_{(r)} $  will be  called the  
{\it $r$-fold 
Massey product of $\xi$ and $a$} and denoted by 
$$
\Big\langle \ \underset{r\ times }{\underbrace {\xi, \ldots, \xi}},\ a\ \Big\rangle
\in MZ_{(r)}\Big/MB_{(r)}.
$$
We say that the {\it length } of this product is equal to $r$.
\end{defi}

\bere\lb{r:defs}
Observe that the indeterminacy of this Massey product is 
less than the indeterminacy of a general Massey product
as defined for example in \cite{Kr}.
Observe also that usually the length of Massey product is defined as the number of the arguments
inside the brackets, so our notion of length is less by 1 than the standard one.
\enre

We obtain a \ho~
$$
MH_{(r)} \arrr {\Delta_r}  MH_{(r)};
\ \ 
a\mapsto 
\Big\langle \ \underset{r\ times }{\underbrace {\xi, \ldots, \xi}},\ a\ \Big\rangle$$
\bepr\label{p:diffs_in_F}
We have $\Delta_r^2=0$, and 
 the  cohomology group  
 $H^*(MH^*_{(r)}, \Delta_r)$ is isomorphic to 
 $MH^*_{(r+1)}$. 
 \enpr
\bede\lb{d:mas} The groups 
 $MH^*_{(r)}$ form therefore a spectral sequence,
 which will be called the {\it Massey spectral sequence}
 associated with $(X, \xi, \l)$,
 and denoted by $\MMMM^*_r(X,\xi,\l)$,
 or just $\MMMM^*_r$ for brevity.
 \end{defi}
 
 \subsection{Formal deformation spectral sequence}
 \label{su:form-def-ss}
 
Another spectral sequence is associated to a formal deformation
 of the twisted DeRham complex $\tilde\Omega^*(X)$.
 Let $\L=\cc[[t]]$, consider a representation
 $\wh\g:\pi_1(X)\to\big(\cc[[t]]\big)^*$
 defined by the formula
 $\wh\g(g)=e^{\langle (\a+t)\xi,g\rangle}.$
 
 Denote by $C^*(X,\wh\g)$ the complex of $\wh\g$-equivariant 
 chains on $\wi X$, and by 
 $H^*(X,\wh\g)$ its cohomology. 
 The exact sequence
 $$0\to C^*(X,\wh\g)\arrr t C^*(X,\wh\g) \to C^*(x, \rho_\l) \to 0$$
 gives rise to an exact couple
\begin{equation}\lb{f:e-c-0}
\xymatrix{
H_*(X, \wh\g) \ar[rr]^t &  & H_*(X, \wh\g) \ar[dl]\\
& H^*(X,\r_\l) \ar[ul]  & \\
} 
\end{equation}
\bede\lb{d:def.sp.s}
The associated spectral sequence will be denoted
by $\EEEE_r^*$ and called 
{\it the formal deformation spectral sequence.}.
\end{defi}

\beth[\cite{KP}, th. 3.6]
The spectral sequences 
$\MMMM^*_r$ and $\EEEE^*_r$
are isomorphic.
\enth
\subsection{Formality and strong formality}
\label{su:form-strongform}

Recall the classical notion of {\it formality } of manifolds,
introduced by D. Sullivan.
Let $X$ be a $\smo$ manifold,
$\Omega^*(X)$ be the differential graded algebra of differential forms on $X$, 
and $\MMM^*(X) $ be a minimal model for $\Omega^*(X)$.
The manifold  $X$ is called {\it formal}
if there is a DGA-\ho~ $\MMM^*(X)\to H^*(X)$ inducing an isomorphism
in cohomology. Compact  K\"ahler manifolds are formal as proved in \cite{DGMS}.

\beth\lb{t:form-degen}(\cite{KP},{\rm\  Th. 3.14})
Assume that $X$ is a  formal \ma.
Then for every $\xi\in H^1(X, \zz)$ 
 the spectral sequences $\EEEE^*_r(X,\xi, 1)$ and $\MMMM_r^*(X,\xi, 1)$
degenerate in their second term. Therefore 
all the Massey products $\langle \xi, \ldots, \xi, x\rangle$
of length $\geq 2$, where $x\in H^*(X, \cc)$ are equal to zero.
\enth

Assume that $X$ is connected and denote $\pi_1(X)$ by $G$.
Denote by $Ch(G)$ the set of all \ho s $G\to \cc^*$.
For a character $\r\in Ch(G)$ 
denote by $E_\r$ the corresponding flat vector bundle over $X$.
Put
$$
\bar \O^*(X) = \underset{\r\in Ch(G)}\bigoplus \O^*(X,E_\r).
$$
The pairing $E_{\r}\otimes E_{\eta}\approx E_{\r\eta}$
induces a natural structure of a differential graded algebra
on the vector space $\bar \O^*(X)$.
\begin{defi}\label{d:hyp_form}
A $\smo$ manifold $X$ is {\it strongly formal} if the differential graded algebra 
$\bar \O^*(X)$ is formal.
\end{defi}

This notion was introduced in the paper of T. Kohno
and the author \cite{KP}, although it  was implicit already in the paper 
of H. Kasuya \cite{Ka}. 
It turns out that K\"ahler compact manifolds are strongly formal
(\cite{KP}, Th. 6.5). This is a  strengthening of the formality 
condition, since  there exist formal, but not strongly formal spaces 
(see H. Kasuya's paper \cite{Ka}, \S 9, Example 1  and also  \cite{KP}, Remark 6.6).

\beth\lb{t:str-form-degen}(\cite{KP},{\rm\  Th. 6.3})
Assume that $X$ is a strongly formal \ma.
Then for every $\xi\in H^1(X, \zz)$ 
and $\l\in \cc^*$ the spectral sequences $\EEEE^*_r(X,\xi, \l)$ and $\MMMM_r^*(X,\xi, \l)$
degenerate in their second term. Therefore 
all the Massey products $\langle \xi, \ldots, \xi, x\rangle$
of length $\geq 2$, where $x\in H^*(X, \rho_\l)$ are equal to zero.
\enth

 \subsection{Terminological conventions}
 \bede\lb{d:stops}
 If a spectral sequence $\CCCC^*_r$ 
 degenerates at sheet $m$ in degree $k$ 
 we write 
 $\s_k(\CCCC_r)=m$.
 \end{defi}
 
 In particular, the number $\s_k(\MMMM(X,\xi,\l))-1$ equals 
 the maximal length of a non-zero Massey product of the form
 $\langle \xi,\ldots , \xi, a\rangle$ where $a\in H^k(X,\rho_\l)$.

 \bede\lb{d:invar-subsp}
 Let $A:L\to L$ be a linear map of a finite-dimensional vector space 
 $L$ over $\cc$, and $\l\in \cc$.
 The $A$-invariant subspace of $L$ corresponding to $\l$
 will be denoted by $N(A,\l)$. The degree of nilpotency of 
 $(A-\l) ~|~ N(A,\l)$ will be denoted by $\nu(A,\l)$. 
 If $N(A,\l)=0$, we set $\nu(A,\l)=0$ by convention.
 Thus the  number $\nu(A,\l)$
 is equal to  the maximal size of a Jordan block of $A$ with eigenvalue $\l$.
\end{defi}

\section{Main theorem}
\lb{s:m_th}
Let $\phi: M\to M$ be a diffeomorphism of a $\smo$ compact 
connected 
\ma,
denote by $X$ its mapping torus. There is a natural fibration
$p:X\to S^1$, denote by $\xi\in H^1(X, \zz)$
the corresponding cohomology class. Let $\l\in \cc^*$.
Composing the  homomorphism $\zz\to \cc^*: n\mapsto e^{n\l}$
with the \ho~ $p_*:\pi_1(X)\to\pi_1(S^1)\approx \zz$
we obtain a representation $\r_\l:\pi_1(X)\to \cc^*$.
Denote by $H^*(X,\r_\l) $ the corresponding
cohomology with local coefficients.
We have a natural pairing
$$ H^*(X,\cc)  \otimes  H^*(X,\r_\l)\to H^*(X,\r_\l) $$
and the corresponding Massey products
$\langle \xi, \ldots \xi, x\rangle$
are defined for $x\in H^*(X,\r_\l)$
(see Section \ref{su:mas-ss}).
Denote by $\m_k(\l)$ the maximal number $r$
such that $\langle \ \underset{r}{\underbrace {\xi, \ldots, \xi}},\ x\ \rangle$
is not equal to zero. 
Then $\m_k(\l)+1$ is the number of the sheet where 
the two spectral sequences discussed at the previous section
degenerate.
Consider the endomorphism 
$\phi^*_k : H^k(M, \cc) \to H^k(M, \cc)$,
and denote by $J_k(\l)$ 
the maximal size of its Jordan block of eigenvalue $\l$.

\beth
\lb{t:main}
\been\item We have $J_k(\l)=\m_k(\l)$ for every $k$ and $\l$.
\item
If $X$ is a strongly formal space (e.g. a compact K\"ahler manifold)
then all Jordan blocks of 
$\phi_k^*$
are of size 1. 
\item 
If $X$ is a formal space, then all Jordan blocks of eigenvalue 1 
are of size 1.
\enen 
\enth
\Prf
The parts 2) and 3) follow immediately from the part 1)
in view of degeneracy of corresponding spectral sequences. 
Proceeding to the part 1), choose any $\a\in \cc$ \sut~ $e^\a=\l$ and 
consider the exact couple $\EEEE$ (see \rrf{f:e-c-0}).

Replacing $t$ by $e^{t+\a}-e^\alpha$ in its upper line we obtain an exact couple $\EEEE'$
isomorphic to $\EEEE$
(since the element $\l(e^t -1)/t\in \cc[[t]]$ is invertible).
Let $L=\cc[u, u^{-1}]\approx \cc[\zz]$, we have a tautological representation
$\zz\to L^*$, composing it with the \ho~ $p_*:\pi_1(X)\to \zz$
we obtain a representation $\b:\pi_1(X)\to L^*$.
The corresponding twisted cohomology of $X$ will be denoted
by 
$H^*(X,\b)$.
We have an exact couple

\begin{equation}\lb{f:e-c-22}
\xymatrix{
H^*(X,\b) \ar[rr]^{u-\l} &  & H^*(X,\b) \ar[dl]\\
& H^*(X,\r_\l) \ar[ul]  & \\
} 
\end{equation}
denote it by $(\DDDD)$.
Consider a ring \ho~
$h:L\to\L; \ \  h(u)=e^{t+\alpha}$.
We have $\wh\g=h\circ \b$, therefore the exact couples 
$\EEEE$ and $\DDDD$ are related by a \ho~ that equals to
$\id$ on the group $H^*(X,\r_\l)$ and to the \ho~ $h$ on the modules 
$H^*(X,\b)$.
This \ho~ induces the identity isomorphism 
on the first term of the corresponding spectral sequences,
therefore they are isomorphic. 
The advantage of the \ses~ $\DDDD$ 
is that it is computable in terms of 
the monodromy \ho~ with the help of the next Lemma.

\bele\lb{l:X_M}
There is a commutative diagram of $L$-\ho s
\begin{equation}\lb{f:x-m}
\xymatrix{
H^k(X,\b) \ar[rr]^{u-\l} &  & H^k(X,\b) \\
 H_{k-1}(M,\cc) \ar[u]^\approx\ar[rr]^{\phi_*-\l}  &  &  H_{k-1}(M,\cc)  \ar[u]^\approx  & \\
} 
\end{equation}
\enle
\Prf
Let $\ove X$ be the corresponding infinite cyclic covering of $X$.
The simplicial chain complex $C_*(\ove{X})$
is a free finitely generated chain complex of $L$-modules.
Its homology $H_*(\ove{X})$
is isomorphic to $H_*(M)$ as $L$-module
(the element $u\in L$ acts as $\phi_*$ on $H_*(M)$).
The universal coefficient theorem implies
an $L$-module isomorphism
$$
H^k(X,\b)
\approx \Ext^1_L( H_{k-1}(M,\cc), L) \approx  H_{k-1}(M,\cc).$$
The exact couple $\DDDD$
is therefore isomorphic to the following one

\begin{equation}\lb{f:e-c-3}
\xymatrix{
H_*(M, \cc) \ar[rr]^{\phi_*-\l} &  & H_*(M, \cc) \ar[dl]^j\\
& H^*(X,\r_\l) \ar[ul]^l  & \\
} 
\end{equation}
Here $\deg i=\deg l=0, \ \deg j=1$.
Let $A_k=N(\phi_*^k, \l)$;
denote by $B_k$ the sum of all subspaces $N(\phi_*^k, \m)$ with
$\mu\not=\l$.
The restriction $(\phi_*^k-\l)~|~A_k$
is nilpotent of degree equal to $J_k(\l)$,
and the restriction $(\phi_*^k-\l)~|~B_k$
is an isomorphism of $B_k$ onto itself.
The assertion of the theorem follows now 
from the following lemma.
\bele\lb{l:decomp}
Let $\EEEE $ be a graded  exact couple:
\begin{equation}\lb{f:e-c-4}
 \xymatrix{
D \ar[rr]^i &  & D \ar[dl]^j\\
& E \ar[ul]^l  & \\
} 
\end{equation}
with $\deg i =\deg l = 0, \ \deg j = 1$.
Assume that  the $k$-th component  $i_{k}:D_{k}\to D_{k}$
of the \ho~ $i$  decomposes as 
$$
\d\oplus\tau: A\oplus B\to  A\oplus B
$$
where $\d$ is nilpotent of degree $m$ and 
$\tau $ is injective.
Then the spectral sequence $\EEEE_*$ degenerates 
at the step $m+1$ in degree $k$.
\enle
\Prf
The proof is an easy diagram chasing.
At the step $r$ the upper line of the derived exact couple $\EEEE_r$
equals $i: i^{r-1}(D)\to i^{r-1}(D)$. 
Observe that $i:i^{m}(D)\to i^{m}(D)$ is injective,
therefore $d_{m+1}:E_{m+1}^k\to E_{m+1}^k$ equals $0$.
On the other hand the \ho~ $i :i^{m-1}(A)\to i^{m-1}(A)$ 
equals $0$, and $i^{m-1}(A)\not=0$, therefore 
the differential $d_m$ in the module $E_{m}^k$ 
is non-zero. $\qs$

\section{Examples}
\label{s:examples}

\subsection{The Heisenberg group}
\label{su:heis}

A classical example of a non-formal space arises from the Heisenberg 
group (see \cite{DGMS}, p. 261). Let $N$ denote the group 
of all upper triangular matrices
$$
\left( 
\begin{matrix}
 1 & a & b \\
  1 & 1 & c \\
   0 & 0 & 1 \\
\end{matrix}
\right)
$$
with real coefficients; let $\G$ be the subgroup 
of matrices with integer coefficients. 
The space $W=N/\G$ is a compact three-dimensional manifold 
and 
\begin{equation}\label{f:heis-betti}
 b_1(W)=b_2(W)=2.
\end{equation}
The group $H^1(W,\rr)$ is generated by elements $ x, y$, and 
$\langle x,x,y\rangle \not=0$.

The space $W$ is fibred over $\ttt^2$ with fiber $S^1$.
Composing this fibration with the projection $\ttt^2\to S^1$ 
we obtain a fibration $W\to S^1$ with fiber $\ttt^2$.
Thus $W$ is a mapping torus with a monodromy $\phi:\ttt^2\to \ttt^2$.
Let $\l\in \cc^*$. The Milnor exact sequence 
$$
\ldots \to  H_i(\ttt^2)\arrr {\phi_*-1}  H_i(\ttt^2) \to H_i(W) \to H_{i-1}(\ttt^2)\to \ldots
$$
together with \rrf{f:heis-betti}
shows immediately that $\phi^1$ has only one eigenvalue of multiplicity 2, namely $1$.
Moreover, $\phi^1$ must have a Jordan block of size 2,
and this corresponds to the above non-vanishing Massey product of length 2. 

\subsection{A formal but not strongly formal mapping torus}
\lb{su:form-tor}

We will now use the main theorem to construct a 
formal but not strongly formal mapping torus.
Let $A$ be an endomorphism of $\zz^n\oplus\zz^n$ defined 
as follows: $A(x,y)=(-x+y, -y)$.
This endomorphism
conserves the standard symplectic pairing 
on $\zz^{2n}$, therefore there is a diffeomorphism
$\phi:M_{2g}\to M_{2g}$ inducing the \ho~ $A$ in $H_1(M_{2g}, \zz)$.
Let $X$ be the mapping torus of $\phi$,
and $p:\ove{X}\to X$ be the infinite cyclic covering.
We have $\ove{X} \sim M_{2g}$, and the Milnor exact sequence
for the covering $p$
reduces to the following:
$$
H_i(M_{2g})\arrr {\phi_*-1}  H_i(M_{2g}) \to H_i(X) \to H_{i-1}(M_{2g}).
$$
It follows that $H^*(X, \cc)\approx H^*(S^1\times S^2, \cc)$
and it is easy to prove that  this space is formal.
However the eigenvalue $-1$ of the map 
$\phi^*:H_i(M_{2g}) \to H_i(M_{2g})$
has a Jordan block of size $n$, therefore it is not a strongly formal space.

\section{A generalization}
\lb{s:gener}

The main theorem of the paper is readily carried over  to a 
somewhat more general (but less geometrically appealing)
case of arbitrary compact manifolds $X$ endowed with
a non-zero cohomology class $\xi\in H^1(X,\zz)$.
Let $p:\ove{X}\to X$ be the corresponding infinite cyclic covering;
$H^k(X, L)$ is a finitely gernerated $L$  module; denote by $T_k$ its $L$-torsion 
submodule.
Then $T_k$ is a finite dimensional vector space over $\cc$.
(In the case when $X$ is a mapping torus of a map $\phi:M\to M$,
we have $T_k\approx H^k(M,\cc)$.)
Denote by $f_k$ the automorphism $T_k\to T_k$ induced 
by the generator $u$ of the structure group of the covering $p$.

\beth\label{t:gener}
Let $\l\in \cc^*$. 

\been\item We have $\nu(f_k, \l)=\mu_k(\l)$.
\item If $X$ is a strongly formal space (e.g. a compact K\"ahler manifold)
then all the Jordan blocks of $f_k$ are of size 1.
\item 
If $X$ is a formal space, then all the Jordan blocks of eigenvalue 1 
of $f_k$ are of size 1.
\enen 
\enth
\Prf
The proof goes on the same lines as the proof of the main theorem.
We will briefly indicate the necessary modifications.
We have an exact couple

\begin{equation}\lb{f:e-c-2}
\xymatrix{
H^*(X,\b) \ar[rr] &  & H^*(X,\b) \ar[dl]\\
& H^*(X,\r_\l) \ar[ul]  & \\
} 
\end{equation}
The $L$-module $H^*(X,\b)$
is isomorphic to a direct sum 
$F_0\oplus F_1\oplus F_2$ where $F_2$ is a free finitely generated 
$L$-module, $F_1$ is a direct sum of cyclic modules
of the form $L/PL$ with $(u-\l)~\nmid~ P$, and $F_0$ is a direct
sum of cyclic modules of the form $L/(u-\l)^rL$.
Put $D_0=F_0, \ D_1=F_1\oplus F_2$.
Contrarily to the case of mapping tori, the $L$-module 
$D_1$ is not a finite-dimensional vector space,
and the map $(u-\l)~|~ D_1: D_1\to D_1$ is not an isomorphism.
However, it is injective, and as before an application of Lemma 
\ref{l:decomp}
completes the proof. $\qs$

\section{Acknowledgements}
\lb{s:acknow}

This work was  finalized during 
author's stay at the Laboratory of Algebraic Geometry 
of HSE (Moscow) in August 2016. The author is  grateful  to 
Professor Fedor Bogomolov for his help and support.

\end{document}